\documentclass[11pt]{article}
\usepackage{amsmath,amssymb,amstext,dsfont,fancyvrb,float,fontenc,graphicx,subfigure, theorem}

\title{\Large\bf On Neumann and Poincare problems\\
for Laplace equation}

\author{\sc Vladimir Ryazanov}
\date{}

\cleardoublepage 
\usepackage[strict]{changepage}
\def\abstractname{Abstract -}   % <-----------------
\def\abstract{\begin{adjustwidth}{1cm}{1cm} \par    \footnotesize \noindent {\bf \abstractname}
\def\endabstract{ \end{adjustwidth} \smallskip }}
{\theorembodyfont{\itshape}}
{\theorembodyfont{\itshape}}
{\theorembodyfont{\itshape}}
{\theorembodyfont{\itshape}}
{\theorembodyfont{\itshape}}
{\theorembodyfont{\rm}}
{\theorembodyfont{\rm}}
{\theorembodyfont{\rm}}
{\theorembodyfont{\rm }}

\begin{document}
\maketitle
\vskip 1.5em

\vskip 1.5em
 \begin{abstract}
It is proved the existence of nonclassical solutions of the Neumann
problem for the harmonic functions in the Jordan rectifiable domains
with arbitrary measurable boundary distributions of normal
derivatives. The same is stated for a special case of the Poincare
problem on directional derivatives. Moreover, it is shown that the
spaces of the found solutions have the infinite dimension.
 \end{abstract}


\begin{MSC}
primary   31A05, 31A20, 31A25, 31B25, 35Q15; se\-con\-da\-ry 30E25,
31C05, 34M50, 35F45.
\end{MSC}

\bigskip

\section{Introduction}

It is well--known that the Neumann problem has no classical
solutions generally speaking even for some continuous boundary data.
The main goal of this short note is to show that the problem has
nonclassical solutions for arbitrary measurable data. The result is
based on a reduction of this problem to the Hilbert
(Riemann-Hilbert) boundary value problem recently solved for
arbitrary measurable coefficients and for arbitrary measurable
boundary data in \cite{R1}, see also \cite{R2} on the Dirichlet
problem.

\medskip

Let us start from a more general {\bf problem on directional
derivatives} in the unit disk $\mathbb D = \{ z\in\mathbb{C}:
|z|<1\}$, $z=x+iy$. The classic setting of the latter problem is to
find a function $u:\mathbb D\to\mathbb R$ that is twice continuously
differentiable, admits a continuous extension to the boundary of
$\mathbb D$ together with its first partial derivatives, satisfies
the Laplace equation
\begin{equation}\label{eqLAPLACE}
\Delta u\ :=\ \frac{\partial^2 u}{\partial x^2}\ +\ \frac{\partial^2
u}{\partial y^2}\ =\ 0 \quad\quad\quad\forall\ z\in\mathbb D
\end{equation} and the boundary condition with a prescribed
continuous date $\varphi : \partial\mathbb D\to\mathbb R$:
\begin{equation}\label{eqDIRECT}
\frac{\partial u}{\partial \nu}\ =\ \varphi(\zeta) \quad\quad\quad
\forall\ \zeta\in\partial\mathbb D
\end{equation}
where $\frac{\partial u}{\partial \nu}$ denotes the derivative of
$u$ at $\zeta$ in a direction $\nu = \nu(\zeta)$, $|\nu(\zeta)|=1$:
\begin{equation}\label{eqDERIVATIVE}
\frac{\partial u}{\partial \nu}\ :=\ \lim_{t\to 0}\
\frac{u(\zeta+t\cdot\nu)-u(\zeta)}{t}\ .
\end{equation}

{\bf The Neumann problem} is a special case of the above problem on
directional derivatives with  the boundary condition
\begin{equation}\label{eqNEUMANN}
\frac{\partial u}{\partial n}\ =\ \varphi(\zeta) \quad\quad\quad
\forall\ \zeta\in\partial\mathbb D
\end{equation}
where $n$ denotes the unit interior normal to $\partial\mathbb D$ at
the point $\zeta$.

\medskip

In turn, the above problem on directional derivatives is a special
case of {\bf the Poincare problem} with  the boundary condition
\begin{equation}\label{eqPOINCARE}
a\cdot u\ +\ b\cdot\frac{\partial u}{\partial \nu}\ =\
\varphi(\zeta) \quad\quad\quad \forall\ \zeta\in\partial\mathbb D
\end{equation}
where $a=a(\zeta)$ and $b=b(\zeta)$ are real-valued functions given
on $\partial\mathbb D$.

\medskip

Recall also that twice continuously differentiable solutions of the
Laplace equation are called {\bf harmonic functions}. As well known,
such functions are infinitely differentiable.

\section{On nonclassical solutions of boundary problems}

Let us start from the more general problem on directional
derivatives.

\medskip

{\bf Theorem 1.} {\it\, Let $\nu:\partial\mathbb D\to\mathbb C$,
$|\nu (\zeta)|\equiv 1$, and $\varphi:\partial\mathbb D\to\mathbb R$
be measurable functions. Then there exist harmonic functions
$u:\mathbb D\to\mathbb R$ such that
\begin{equation}\label{eqLIMIT} \lim\limits_{z\to\zeta}\ \frac{\partial u}{\partial \nu}\ (z)\ =\
\varphi(\zeta)\end{equation} along any nontangential paths to a.e.
point $\zeta\in\partial\mathbb D$.}

\bigskip

{\bf Remark 1.} We are able to say more in the case of $\mathrm
{Re}\ n(\zeta)\cdot\overline{\nu(\zeta)}>0$. Indeed, the latter
magnitude is a scalar product of $n=n(\zeta)$ and $\nu =\nu(\zeta)$
interpreted as vectors in $\mathbb R^2$ and it has the geometric
sense of projection of the vector $\nu$ onto the inner normal $n$ to
$\partial\mathbb D$ at the point $\zeta$. In view of
(\ref{eqLIMIT}), since the limit $\varphi(\zeta)$ is finite, there
is a finite limit $u(\zeta)$ of $u(z)$ as $z\to\zeta$ in $\mathbb D$
along the straight line passing through the point $\zeta$ and being
parallel to the vector $\nu$ because along this line
\begin{equation}\label{eqDIFFERENCE} u(z)\ =\ u(z_0)\ -\ \int\limits_{0}\limits^{1}\
\frac{\partial u}{\partial \nu}\ (z_0+\tau (z-z_0))\ d\tau\
.\end{equation} Thus, at each point with condition (\ref{eqLIMIT}),
there is the directional derivative
\begin{equation}\label{eqPOSITIVE}
\frac{\partial u}{\partial \nu}\ (\zeta)\ :=\ \lim_{t\to 0}\
\frac{u(\zeta+t\cdot\nu)-u(\zeta)}{t}\ =\ \varphi(\zeta)\ .
\end{equation}

\bigskip

In particular, in the case of the Neumann problem, $\mathrm {Re}\
n(\zeta)\cdot\overline{\nu(\zeta)}\equiv 1>0$ and we have by Theorem
1 and Remark 1 the following significant result.

\medskip

{\bf Theorem 2.} {\it\, For each measurable function
$\varphi:\partial\mathbb D\to\mathbb R$, one can find harmonic
functions $u:\mathbb D\to\mathbb R$ such that, at a.e. point
$\zeta\in\partial\mathbb D$, there exist:

\bigskip

1) the finite radial limit
\begin{equation}\label{eqLIMIT1}
u(\zeta)\ :=\ \lim\limits_{r\to 1}\ u(r\zeta)\end{equation}

2) the normal derivative
\begin{equation}\label{eqNORMAL}
\frac{\partial u}{\partial n}\ (\zeta)\ :=\ \lim_{t\to 0}\
\frac{u(\zeta+t\cdot n)-u(\zeta)}{t}\ =\ \varphi(\zeta)
\end{equation}

3) the nontangential limit
\begin{equation}\label{eqLIMIT2} \lim\limits_{z\to\zeta}\ \frac{\partial u}{\partial n}\ (z)\ =\
\frac{\partial u}{\partial n}\ (\zeta)\end{equation} where
$n=n(\zeta)$ denotes the unit interior normal to $\partial\mathbb D$
at the point $\zeta$.}

\bigskip

\begin{proof} To prove Theorem 1, let us show that the problem on
directional derivatives is equivalent to the corresponding
Riemann-Hilbert problem.

\medskip

Indeed, let $u$ be a harmonic function $u:\mathbb D\to\mathbb R$
satisfying the boundary condition (\ref{eqLIMIT}). Then the
functions $U=u_x$ and $V=-u_y$ satisfy the system of Cauchy-Riemann:
$U_y=-V_x$ and $U_x=V_y$ in view of (\ref{eqLAPLACE}). Thus, the
function $f=U+iV$ is analytic in $\mathbb D$ and along any
nontangential path to a.e. $\zeta\in\partial\mathbb D$
\begin{equation}\label{eqRH} \lim\limits_{z\to\zeta}\ \mathrm {Re}\
\nu(\zeta)\cdot f(z)\ =\ \varphi(\zeta)\end{equation} that is
equivalent to (\ref{eqLIMIT}). Inversely, let $f:\mathbb D\to\mathbb
C$ be an analytic function satisfying the boundary condition
(\ref{eqRH}). Then any indefinite integral $F$ of $f$ is also a
single-valued analytic function in $\mathbb D$ and $u=\mathrm {Re}\
F$ is a harmonic function satisfying the boundary condition
(\ref{eqLIMIT}) because the directional derivative
\begin{equation}\label{eqBOUNDARY} \frac{\partial u}{\partial \nu}\
=\ \mathrm {Re}\ \overline{\nu}\cdot\nabla u\ =\ \mathrm {Re}\
\nu\cdot\overline{\nabla u}\ =\ (\nu,\nabla u)\end{equation} is the
scalar product of $\nu$ and the gradient $\nabla u$ interpreted as
vectors in $\mathbb R^2$.

\medskip

Thus, Theorem 1 is a direct consequence of Theorem 2.1 in \cite{R1}
on the Riemann-Hilbert problem with
$\lambda(\zeta)=\overline{\nu(\zeta)}$, $\zeta\in\partial\mathbb D$.
\end{proof}

\bigskip

The following result in domains bounded by rectifiable Jordan curves
is proved perfectly similar to Theorem 1 but it is based on more
general Theorem 3.1 in \cite{R1}.

\medskip

{\bf Theorem 3.} {\it\, Let $D$ be a domain in $\mathbb C$ bounded
by a rectifiable Jordan curve, $\nu:\partial D\to\mathbb C$, $|\nu
(\zeta)|\equiv 1$, and $\varphi:\partial D\to\mathbb R$ be
measurable functions with respect to the natural parameter. Then
there exist harmonic functions $u: D\to\mathbb R$ such that along
any nontangential paths
\begin{equation}\label{eqLIMIT-R} \lim\limits_{z\to\zeta}\ \frac{\partial u}{\partial \nu}\ (z)\ =\
\varphi(\zeta)\end{equation}  for a.e. point $\zeta\in\partial D$
with respect to the natural parameter.}

\bigskip

{\bf Remark 2.} Again we are able to say more in the case with
$\mathrm {Re}\ n\cdot\overline{\nu}>0$ where $n=n(\zeta)$ is the
unit inner normal at a point $\zeta\in\partial D$ with a tangent to
$\partial D$. In view of (\ref{eqLIMIT-R}), since the limit
$\varphi(\zeta)$ is finite, there is a finite limit $u(\zeta)$ of
$u(z)$ as $z\to\zeta$ in $\mathbb D$ along the straight line passing
through the point $\zeta$ and being parallel to the vector $\nu$
because along this line, for $z$ and $z_0$ that are close enough to
$\zeta$,
\begin{equation}\label{eqDIFFERENCE-R} u(z)\ =\ u(z_0)\ -\ \int\limits_{0}\limits^{1}\
\frac{\partial u}{\partial \nu}\ (z_0+\tau (z-z_0))\ d\tau\
.\end{equation} Thus, at each point with the condition
(\ref{eqLIMIT-R}), there is the directional derivative
\begin{equation}\label{eqPOSITIVE-R}
\frac{\partial u}{\partial \nu}\ (\zeta)\ :=\ \lim_{t\to 0}\
\frac{u(\zeta+t\cdot\nu)-u(\zeta)}{t}\ =\ \varphi(\zeta)\ .
\end{equation}

%\bigskip

In particular, in the case of the Neumann problem, $\mathrm {Re}\
n(\zeta)\cdot\overline{\nu(\zeta)}\equiv 1>0$ and we have by Theorem
3 and Remark 2 the following significant result. Here we also apply
the well-known fact that any rectifiable curve has a tangent a.e.
with respect to the natural parameter.

\medskip

{\bf Theorem 4.} {\it\, Let $D$ be a domain in $\mathbb C$ bounded
by a rectifiable Jordan curve and $\varphi:\partial D\to\mathbb R$
be a measurable function with respect to the natural parameter. Then
one can find harmonic functions $u: D\to\mathbb R$ such that, at
a.e. point $\zeta\in\partial D$ with respect to the natural
parameter, there exist:

\bigskip

1) the finite normal limit
\begin{equation}\label{eqLIMIT1-R}
u(\zeta)\ :=\ \lim\limits_{z\to\zeta}\ u(z)\end{equation}

2) the normal derivative
\begin{equation}\label{eqNORMAL-R}
\frac{\partial u}{\partial n}\ (\zeta)\ :=\ \lim_{t\to 0}\
\frac{u(\zeta+t\cdot n)-u(\zeta)}{t}\ =\ \varphi(\zeta)
\end{equation}

3) the nontangential limit
\begin{equation}\label{eqLIMIT2-R} \lim\limits_{z\to\zeta}\ \frac{\partial u}{\partial n}\ (z)\ =\
\frac{\partial u}{\partial n}\ (\zeta)\end{equation} where
$n=n(\zeta)$ denotes the unit interior normal to $\partial D$ at the
point $\zeta$.}

\medskip

Note that here the tangent $\tau(s)$ to $\partial D$ is measurable
with respect to the natural parameter $s$ as the derivative
$d\zeta(s)/ds$ and, thus, the inner normal $n(s)$ to $\partial D$ is
also measurable with respect to the natural parameter.

\bigskip

{\bf Remark 3.} Similarly, on the basis of Theorem 4.1 in \cite{R1},
in the case of arbitrary Jordan domains $D$ in $\mathbb C$, we are
able to conclude that there exist harmonic functions $u:D\to\mathbb
R$ for which the limit relation (\ref{eqLIMIT-R}) holds in the sense
of the unique principal asymptotic value for a.e. $\zeta\in\partial
D$ with respect to the harmonic measure, see the corresponding
definitions and comments in \cite{R1}, Section 4.

\bigskip

%\cc
\section{On the dimension of spaces of solutions}

\bigskip

Finally, we have the following significant result.

\bigskip

{\bf Theorem 5.}{\it\, The spaces of harmonic functions in Theorems
1-4, being nonclassical solutions of the problem on directional
derivatives and the Neumann problem, have the infinite dimension for
any prescribed measurable boundary data.}

\bigskip

\begin{proof} In view of the equivalence of the problem on the
directional derivatives to the corresponding Hilbert
(Riemann-Hilbert) boundary value problem established under the proof
of Theorem 1, the conclusion of Theorem 5 follows directly from
Theorem 5.2 and Remark 5.2 in \cite{R1}.
\end{proof}

\bigskip

% \section*{Acknowledgments} Please write at the end of the paper the acknowledgements.

 %%%%%%%% AUTHORS

 { \footnotesize
\medskip
\medskip
 \vspace*{1mm}

\noindent {\it Vladimir Ryazanov}\\
Institute of Applied Mathematics and Mechanics \\
of National Academy of Sciences of Ukraine, \\
19 gen. Batyuk Str., 84116 Slavyansk, UKRAINE\\
E-mail: {\tt vl$\underline{\ \ }$\,ryazanov@mail.ru}}

 \end{document}